\documentclass[a4paper,11pt]{article}
\usepackage{amssymb}
\usepackage{amsfonts}
\usepackage{amsmath}

\input{tcilatex}
\begin{document}

\begin{titlepage}
\title{\bf  Mechanical Equations on Bi-Para Conformal Geometry}
\author{Zeki Kasap {\footnote{Corresponding Author. E-mail: zkasap@pau.edu.tr; 
Tel:+9 0 258 296 11 82; Fax:+9 0 258 296 12 00}} \\
  {\small Department of Elementary Education, Faculty of Education,  Pamukkale University}
\\{\small  20070 Denizli, Turkey.}
\\ Mehmet Tekkoyun \footnote{Author. E-mail: tekkoyun@pau.edu.tr; Tel: +9 0 258 296 36 16} \\
 {\small Department of Mathematics, Pamukkale University}
\\{\small  20070 Denizli, Turkey.} }
\date{\today}
\maketitle
\begin{abstract}

This study is an extented analogue to conformal geometry of the paper given by [14]. 
    
Also, the geometric and physical results related to bi-para-conformal-dynamical systems are also presented.

{\bf Keywords:} Conformal geometry, bi-para Lagrangian, bi-para Hamiltonian.

{\bf MSC:} 53C15, 70H03, 70H05.
\end{abstract}
\end{titlepage}

\section{Introduction}

There are a great number of applications for differential geometry and
mathematical physics. This applications can be use in many areas in this
century. One of the most important applications of differential geometry is
on geodesics. A geodesic is the shortest route between two points. Geodesics
can be found with the help of the Euler-Lagrange and Hamilton equations.
Also, the information about them can be seen in many mechanical and geometry
books. It is well known that differential geometry provides a suitable field
for studying Lagrangians and Hamiltonians of classical mechanics\ and field
theory. So, the dynamic equations for moving bodies were obtained for
Lagrangian and Hamiltonian mechanics by many authors and are illustrated as
follows:

\textbf{I.} \textbf{Lagrange Dynamics Equations} \cite{klein,deleon,abraham}:%
\textbf{\ }Let $M$ be an $n$-dimensional manifold and $TM$ its tangent
bundle with canonical projection $\tau _{M}:TM\rightarrow M$. $TM$ is called
the phase space of velocities of the base manifold $M$. Let $L:TM\rightarrow
R$ be a differentiable function on $TM$\ and is called the \textbf{%
Lagrangian function}. We consider closed 2-form on $TM$:%
\begin{equation}
\Phi _{L}=-dd_{J}L.  \label{0.1}
\end{equation}%
(if $J^{2}=-I,$ $J$\ is a complex structure or if $J^{2}=I,$ $J$ is a
para-complex structure and $Tr(J)=0$). Consider the equation%
\begin{equation}
\mathbf{i}_{\xi }\Phi _{L}=dE_{L}.  \label{1}
\end{equation}%
Where the semispray $\xi $ is a vector field. We know that $E_{L}=V\left(
L\right) -L$\ \ is an energy function and $V=J(\xi )$\ a Liouville vector
field. Here $dE_{L}$ denotes the differential of $E_{L}$. It is well-known
that (\ref{1}) under a certain condition on $\xi $ is the intrinsical
expression of the Euler-Lagrange equations of motion. This equation is named
\ as \textbf{Lagrange dynamical equation}. The triple $(TM,\Phi _{L},\xi )$
is known as \textbf{Lagrangian system} on the tangent bundle $TM$. The
operations run on (\ref{1}) for\ any coordinate system $(q^{i}(t),p_{i}(t))$%
. Infinite dimension \textbf{Lagrangian's\ equation }is obtained the form
below:%
\begin{equation}
\begin{array}{l}
\frac{dq^{i}}{dt}=\dot{q}^{i}\ ,\ \frac{d}{dt}\left( \frac{\partial L}{%
\partial \dot{q}^{i}}\right) =\frac{\partial L}{\partial q^{i}}\ ,\
i=1,...,n.%
\end{array}
\label{2}
\end{equation}%
\textbf{II.} \textbf{Hamilton Dynamics Equations }\cite{deleon,abraham}%
\textbf{: }Let $M$ be the base manifold and its cotangent manifold $T^{\ast
}M$. By a symplectic form we mean a 2-form $\Phi $ on $T^{\ast }M$ such that:

\textbf{(i)} $\Phi $ is closed , that is, $d\Phi =0;$ \textbf{(ii)} for each 
$z\in T^{\ast }M$ , $\Phi _{z}:T_{z}T^{\ast }M\times T_{z}T^{\ast
}M\rightarrow 
\mathbb{R}
\ $is weakly non$\deg $enerate. If $\Phi _{z}$ in (ii) is non$\deg $enerate,
we speak of a \textbf{strong symplectic form}. If (ii) is dropped we refer
to $\Phi $\ as a presymplectic form\textbf{. }Now\textbf{\ }let $(T^{\ast
}M,\Phi )$ us take as a symplectic manifold. A vector field $Z_{H}:T^{\ast
}M\rightarrow TT^{\ast }M$ is called \textbf{Hamiltonian vector field}\ if
there is a $C^{1}$ \textbf{Hamiltonian function}\ $H:T^{\ast }M\rightarrow 
\mathbb{R}
$ such that \textbf{Hamilton dynamical equation} is determined by%
\begin{equation}
\mathbf{i}_{Z_{H}}\Phi =dH.  \label{3}
\end{equation}%
We say that $Z_{H}$ is locally Hamiltonian vector field if $\Phi $\ is
closed. Where $\Phi $ shows the canonical symplectic form so that $\Phi
=-d\lambda ,$ $\lambda =J^{\ast }(\omega ),$ such that $J^{\ast }$ a dual of 
$J,$ $\omega $ a 1-form on $T^{\ast }M$. The triple $(T^{\ast }M,\Phi
,Z_{H}) $ is named \textbf{Hamiltonian system }which it is defined on the
cotangent bundle $T^{\ast }M.$ From the local expression of $Z_{H}$ we see
that $(q^{i}(t),p_{i}(t))$ is an integral curve of $Z_{H}$ iff \textbf{%
Hamilton's equations }are expressed as follows:%
\begin{equation}
\dot{q}^{i}=\frac{\partial H}{\partial p_{i}}\ ,\ \dot{p}_{i}=-\frac{%
\partial H}{\partial q^{i}}.  \label{4}
\end{equation}%
\textbf{Considering information the above, in a lot of articles and books,
it is possible to show how differential geometric methods are applied in
Lagrangian's and Hamiltonian's mechanics in the below. Some works in
paracomplex geometry are used for mathematical models. }

\textit{Cruceanu}, \textit{Fortuny} and \textit{Geada} have presented
paracomplex geometry which is related to algebra of paracomplex number and
the study of the structures on differentiable manifolds called paracomplex
structures. Furthermore, they have considered a compatible neutral
pseudo-Riemannian metric, the para-Hermitian and para-Kahler structures, and
their variants \textbf{\cite{cruceanu}. }Kaneyuki and Kozai have introduced
a class of affine symmetric spaces, which are called para-Hermitian
symmetric spaces, a paracomplex analogue Hermitian symmetric space \textbf{%
\cite{kaneyuki}. }In the study of para-Kahlerian manifolds, \textit{Tekkoyun}
has introduced paracomplex analogues of the Euler-Lagrange and Hamilton
equations. Furthermore, the geometric results on the related mechanic
systems have been presented\textbf{\ \cite{tekkoyun1}. }\textit{Etayo} and%
\textit{\ Santamaria} studied connections attached to non-integrable almost
biparacomplex manifolds. Manifolds endowed with three foliations pairwise
transversal are called $3$-webs. Similarly, they can be algebraically
defined as biparacomplex or complex product manifolds, i.e., manifolds
endowed with three tensor fields of type $(1,1),$ $F,$ $P$ and $J=FoP$,
where the two first are product and the third one is complex, and they
mutually anti-commute. In this case, it is well known that there exists a
unique torsion-free connection parallelizing the structure. A para-K\"{a}%
hlerian manifold $M$\ is said to be endowed with an almost
bi-para-Lagrangian structure (a bi-para-Lagrangian manifold) if $M$\ has two
transversal Lagrangian distributions (involutive transversal Lagrangian
distributions) $D_{1}$\ and $D_{2}$\ \cite{etayo}. \textit{Carinena} and 
\textit{Ibort} obtained the Lax equations which are associated with a
dynamical endowed with a bi-Lagrangian connection and a closed two-form $%
\Omega $ parallel along dynamical field $\Gamma $. The case of Lagrangian
dynamical systems is analysed and the nonnoether constants of motion found
by Hojman and Harleston are recovered as being associated to a reduced Lax
equation. Completely integrable dynamical systems have been shown to be a
particular case of these systems by their \textbf{\cite{carinena}. }\textit{%
Gordejuela} and \textit{Santamaria} have proved that the canonical
connection of a bi-Lagrangian manifold introduced which was by Hess is a
Levi-Civita connection by showing that a bi-Lagrangian manifold (i.e. a
symplectic manifold endowed with two transversal Lagrangian foliations) is
endowed with a canonical semi-Riemannian metric \textbf{\cite{etayo2}. }%
\textit{Kanai} has been concerned with closed $C^{\infty }$ Riemannian
manifolds of negative curvature whose geodesic flows have $C^{\infty }$
stable and unstable foliations. In particular, we have shown that the
geodesic flow of such a manifold is isomorphic to that of a certain closed
Riemannian manifold of constant negative curvature if the dimension of the
manifold is greater than two and if the sectional curvature lies between $%
\frac{-9}{4}$ and $-1$ strictly \textbf{\cite{kanai}. }Since they have shown
fundamental physical properties in turbulence (conservation laws, wall laws,
Kolmogorov energy spectrum,...), symmetries are used to analyse common
turbulence models. A class of symmetry preserving turbulence models has been
proposed. This class has been refined in such a way that the models respect
the second law of thermodynamics. Moreover, an example of model belonging to
the class has been numerically tested by \textit{Razafindralandy} and 
\textit{Hamdouni} \textbf{\ \cite{dina}. }A base-equation method has been
implemented to actualize the hereditary algebra of the Korteweg--de Vries
(KdV) hierarchy and the N-soliton manifold is reconstructed. The novelty of
our approach is the fact that it can in a rather natural way, predict other
nonlinear evolution equations which admit local conservation laws.
Significantly enough, base functions themselves are found to provide a basis
to regard the KdV-like equations as higher order degenerate bi-Lagrangian
systems by \textit{Chakrabarti} and \textit{Talukdar}\textbf{\ \cite%
{chakrabarti}. }Bi-para-complex analogue of Lagrangian and Hamiltonian
systems has been introduced on Lagrangian distributions by \textit{Tekkoyun}
and \textit{Sari}. Additionally, the geometric and physical results related
to bi-para-dynamical systems have also been presented by them \textbf{\cite%
{tekkoyun2010}. }Authors introduced generalized-quaternionic K\"{a}hler
analogue of Lagrangian and Hamiltonian mechanical systems. Moreover , the
geometrical-physical results which are related to generalized-quaternionic K%
\"{a}hler mechanical systems have also been also given by Tekkoyun and Yayli%
\textbf{\ \cite{tekkoyun2011}.}

In the above studies; although para-complex mechanical systems were analyzed
successfully in relatively broad area of science, they have not dealt with
bi-para-complex conformal mechanical systems on the bi-Lagrangian conformal
manifold. In this study, therefore, equations related to bi-para-conformal
mechanical systems on the bi-Lagrangian conformal manifold used in obtaining
geometric quantization have been presented.

\section{Preliminaries}

In this study, all the manifolds and geometric objects are $C^{\infty }$\
and the Einstein summation convention $\left( \sum x_{i}=x_{i}\right) $ is
in use. Also, $A$, $F(TM)$, $\chi (TM)$\ and $\Lambda ^{1}(TM)$\ denote the
set of para-complex numbers, the set of para-complex functions on $TM$, the
set of para-complex vector fields on $TM$\ and the set of para-complex
1-forms on $TM$, respectively. The definitions and geometric structures on
the differential manifold $M$\ given in \cite{cruceanu} may be extended to $%
TM$\ as follows:

\section{Conformal Geometry}

A conformal map is a function which preserves angles. Conformal maps can be
defined between domains in higher dimensional Euclidean spaces, and more
generally on a Riemann or semi-Riemann manifold. A conformal manifold is a
differentiable manifold equipped with an equivalence class of (pseudo)
Riemann metric tensors, in which two metrics $g^{\prime }$ and $g$ are
equivalent if and only if%
\begin{equation}
g^{\prime }=\lambda ^{2}g  \label{4.1}
\end{equation}%
where $\lambda >0$ is a smooth positive function. An equivalence class of
such metrics is known as a conformal metric or conformal class \cite{wiki}.
Two Riemann metrics $g$ and $g^{\prime }$ on $M$ \ are said to be equivalent
if and only if%
\begin{equation}
g^{\prime }=e^{\lambda }g  \label{7}
\end{equation}%
where $\lambda $ is a smooth function on $M$. The equation given by (\ref{7}%
) is called a \textbf{Conformal Structure }\cite{folland}.

\section{Bi-Para-Complex Geometry}

Let $M$ be a differentiable manifold. An almost bi-para-complex structure on 
$M$ is denoted by two tensor fields $F$ and $P$ of type (1,1) giving $%
F^{2}=P^{2}=1$, $F\circ P+P\circ F=0$ \cite{etayo}. It is seen that $P\circ
F $ is an almost complex structure. If the matrix-structure defined by the
almost bi-para-complex structure is integrable then for every point $p\in M$
there exists an open neighbourhood $U$ of $p$ and local coordinates $%
(U;x^{_{i}},y^{i})$ such that%
\begin{eqnarray}
F(\partial /\partial x^{i}) &=&\partial /\partial y^{i},F(\partial /\partial
y^{i})=\partial /\partial x^{i},  \label{2.0} \\
P(\partial /\partial x^{i}) &=&\partial /\partial x^{i},\text{ }P(\partial
/\partial y^{i})=-\partial /\partial y^{i},\text{ }\forall i=\overline{1,n} 
\notag
\end{eqnarray}%
\cite{etayo1}. The existence of these kind of local coordinates on $M$
permits to construct holomorphic local coordinates, $(U;z^{k}),$ $%
z^{k}=x^{k}+\mathbf{i}y^{k}$, $\mathbf{i}^{2}=-1,$ $k=\overline{1,n},$ or
para-holomorphic local coordinates, $(U;z^{k}),$ $z^{k}=x^{k}+\mathbf{j}%
y^{k},k=\overline{1,n},$ $\mathbf{j}^{2}=1$ \cite{gadea, newlander}. $%
(M,g,J) $ is a para-K\"{a}hlerian manifold that always has two transversal
distributions defined by the eigen-spaces associated to $+1$ and $-1$
eigenvalues of $J$. Besides, the mentioned distributions are involutive
Lagrangian distributions if somebody thinks of the symplectic form $\Phi $
defined by 
\begin{equation*}
\Phi (X,Y)=g(JX,Y),\forall X,Y\in \chi (M).
\end{equation*}%
Consider $(x^{i},y^{i})$ to be a real coordinate system on a neighborhood $U$
of any point $p$ of $M.$ Also let $\{\frac{\partial }{\partial x^{i}},\frac{%
\partial }{\partial y^{i}}\}$ and $\{dx^{i},dy^{i}\}$ be natural bases over $%
R$ of the tangent space $T_{p}(M)$ and the cotangent space $T_{p}^{\ast }(M)$
of $M,$ respectively. Then the below equalities may be written by%
\begin{equation}
J(\frac{\partial }{\partial x^{i}})=\frac{\partial }{\partial y^{i}}%
,\,\,\,\,\,\,\,J(\frac{\partial }{\partial y^{i}})=\frac{\partial }{\partial
x^{i}}.  \label{2.1}
\end{equation}%
Let $\ z^{i}=\ x^{i}+\mathbf{j}\ y^{i},$ $\mathbf{j}^{2}=1,$ also be a
para-complex local coordinate system on $M.$ So the vector fields will be
shown:%
\begin{equation}
\frac{\partial }{\partial z^{i}}=\frac{1}{2}\{\frac{\partial }{\partial x^{i}%
}-\mathbf{j}\frac{\partial }{\partial y^{i}}\},\,\,\,\frac{\partial }{%
\partial \overline{z}^{i}}=\frac{1}{2}\{\frac{\partial }{\partial x^{i}}+%
\mathbf{j}\frac{\partial }{\partial y^{i}}\}.  \label{2.2}
\end{equation}%
which represent the bases of $M$. Also, the dual covector fields are 
\begin{equation}
dz^{i}=dx^{i}+\mathbf{j}dy^{i},\,\,\,\,\,d\overline{z}^{i}=dx^{i}-\mathbf{j}%
dy^{i}  \label{2.22}
\end{equation}%
which represent the cobases of $M$. Then the following expression can be
written%
\begin{equation}
J(\frac{\partial }{\partial z^{i}})=\mathbf{-j}\frac{\partial }{\partial 
\overline{z}^{i}},\,\,\,\,\,\,J(\frac{\partial }{\partial \overline{z}^{i}})=%
\mathbf{j}\frac{\partial }{\partial z^{i}}.  \label{2.3}
\end{equation}%
The dual endomorphism $J^{\ast }$ of $T_{p}^{\ast }(M)$ at any point $p$ of
the manifold $M$ satisfies that $J^{\ast 2}=I,$ and is denoted by%
\begin{equation}
J^{\ast }(dz^{i})=\mathbf{-j}d\overline{z}^{i},\,\,\,\,\,\,J^{\ast }(d%
\overline{z}^{i})=\mathbf{j}dz^{i}.  \label{2.4}
\end{equation}%
Let $V^{A}$ be a commutative group $(V,+)$ endowed with a structure of
unitary module over the ring $A.$ Let $V^{R}$ denote the group $(V,+)$
endowed with the structure of real vector space inherited from the
restriction of scalars to $R\mathbf{.}$ The vector space $V^{R}$ will then
be called the real vector space associated to $V^{A}.$ Setting%
\begin{equation}
J(u)=ju,\,\,\,\,\,\,P^{+}(u)=e^{+}u,\,\,\,P^{-}(u)=e^{-}u\,\,,\,\,\,u\in
V^{A},  \label{2.5}
\end{equation}%
the equalities%
\begin{equation}
\begin{array}{l}
J^{2}=1_{V},\,\,\,\,P^{+2}=P^{+},\,\,\,\,P^{-2}=P^{-},\,\,\,\,\,\,P^{+}\circ
\,\,\,P^{-}=P^{-}\circ P^{+}=0 \\ 
P^{+}+P^{-}=1_{V},\,\,\,P^{+}-\,\,\,P^{-}=J, \\ 
P^{-}=(1/2)(1_{V}-J),\,\,\,\,\,\,P^{+}=(1/2)(1_{V}+J), \\ 
j^{2}=1,\,\,\,\,e^{+2}=e^{+},\,\,\,\,e^{-2}=e^{-},\,\,\,\,\,\,e^{+}\circ
\,\,\,e^{-}=e^{-}\circ e^{+}=0,\, \\ 
\,\,e^{+}+e^{-}=1,\,\,\,e^{+}-\,\,\,e^{-}=j, \\ 
e^{-}=(1/2)(1-j),\,e^{+}=(1/2)(1+j).%
\end{array}
\label{2.6}
\end{equation}%
can be found. Also, we calculated which%
\begin{equation}
\begin{array}{l}
P^{\mp }\left( \frac{\partial }{\partial z^{i}}\right) =-e^{\mp }\frac{%
\partial }{\partial \overline{z}^{i}}\text{ \ \ \ \ },\text{ \ \ \ \ }P^{\mp
}\left( \frac{\partial }{\partial \overline{z}^{i}}\right) =e^{\mp }\frac{%
\partial }{\partial z^{i}}, \\ 
P^{\ast \mp }\left( dz^{i}\right) =-e^{\mp }d\overline{z}^{i}\text{ \ \ \ \ }%
,\text{ \ \ \ \ }P^{\ast \mp }\left( d\overline{z}^{i}\right) =e^{\mp
}dz^{i}.%
\end{array}
\label{2.8}
\end{equation}%
If the conformal manifold $(M,g,J=P^{+}-P^{-})$ satisfies the following
conditions simultaneously then the conformal manifold is an almost
para-conformal\textbf{\ }Hermitian manifold. The first expression can be
given as follows:%
\begin{equation}
g(X,Y)+g(X,Y)=0\Leftrightarrow g(X,Y)=0,\,\,\,\,\forall X,Y\in \chi (D_{1}),
\label{2.9}
\end{equation}%
since $P^{+}$ and $P^{-}$ are the projections over $D_{1}$ and $D_{2}$
respectively. Then $(P^{+}-P^{-})(X)=P^{+}X-P^{-}X=P^{+}X=X,$ $%
(P^{+}-P^{-})(Y)=P^{+}Y-P^{-}Y=P^{+}Y=Y.$ Similarly the second expression
can be shown as follows:%
\begin{equation}
g(X,Y)+g(X,Y)=0\Leftrightarrow g(X,Y)=0,\,\,\,\,\forall X,Y\in \chi (D_{2}).
\label{2.10}
\end{equation}%
Let $X=X_{1}+X_{2},Y=Y_{1}+Y_{2}$ be vector fields on $M$ such that $%
X_{1},Y_{1}\in D_{1}$ and $X_{2},Y_{2}\in D_{2}.$ Then%
\begin{equation}
\begin{array}{c}
g(JX,Y)=g(JX_{1}+JX_{2},Y_{1}+Y_{2})=g(X_{1}-X_{2},Y_{1}+Y_{2}) \\ 
=g(X_{1},Y_{1})-g(X_{2},Y_{1})+g(X_{1},Y_{2})-g(X_{2},Y_{2}) \\ 
=-g(X_{2},Y_{1})+g(X_{1},Y_{2}), \\ 
g(X,JY)=g(X_{1}+X_{2},JY_{1}+JY_{2})=g(X_{1}+X_{2},Y_{1}-Y_{2}) \\ 
=g(X_{1},Y_{1})+g(X_{2},Y_{1})-g(X_{1},Y_{2})-g(X_{2},Y_{2}) \\ 
=g(X_{2},Y_{1})-g(X_{1},Y_{2}),%
\end{array}
\label{2.101}
\end{equation}%
and hence $%
g(JX,Y)+g(X,JY)=-g(X_{2},Y_{1})+g(X_{1},Y_{2})+g(X_{2},Y_{1})-g(X_{1},Y_{2})=0, 
$ for all vector fields $X,Y$ on $M$. If the conditions (\ref{2.9}) and (\ref%
{2.10}) are true then $D_{1}$ and $D_{2}$ are Lagrangian distributions in
terms of the 2- form $\Phi (X,Y)=g(JX,Y).$ Therefore, if the almost
para-complex structure $J$\ is integrable then $(M,g,J)$ is a para-conformal
K\"{a}hlerian manifold, or equivalently, $(M,\Phi ,D_{1},D_{2})$ is a
bi-Lagrangian conformal manifold. \cite{gilkey,weyl,kim}. Where $W_{\pm }$
is a conformal para-complex structure to be similar to an integrable almost
(para)-complex $P^{\mp }$ given in (\ref{2.8}). Similarly $W_{\pm }^{\ast }$
are the dual of $W_{\pm }$ structures. So, we adapt the following equations
using (\ref{7}):%
\begin{equation}
\begin{array}{l}
W^{\mp }\left( \frac{\partial }{\partial z^{i}}\right) =-e^{\mp }e^{\lambda }%
\frac{\partial }{\partial \overline{z}^{i}}\text{ \ \ \ \ },\text{ \ \ \ \ }%
W^{\mp }\left( \frac{\partial }{\partial \overline{z}^{i}}\right) =e^{\mp
}e^{-\lambda }\frac{\partial }{\partial z^{i}}, \\ 
W^{\ast \mp }\left( dz_{i}\right) =-e^{\mp }e^{\lambda }d\overline{z}_{i}%
\text{\ \ \ \ },\text{ \ \ \ \ }W^{\ast \mp }\left( dz_{i}\right) =e^{\mp
}e^{-\lambda }d\overline{z}_{i}.%
\end{array}
\label{2.11}
\end{equation}

\section{Conformal Bi-Para Euler-Lagrangians}

Here, conformal bi-para-Euler-Lagrange equations and a conformal
bi-para-mechanical system will be obtained under the consideration of the
basis $\{e^{+},e^{-}\}$ on the bi-Lagrangian conformal manifold $(M,\Phi
,D_{1},D_{2})$. Let $(W^{+},W^{-})$ be an almost bi-para-complex conformal
structure on $(M,\Phi ,D_{1},D_{2})$, and $(z^{i},\overline{z}^{i})$ be its
para-complex coordinates. Let the vector field $\xi $ be a semispray given by%
\begin{equation}
\begin{array}{l}
\xi =e^{+}(\xi ^{i+}\frac{\partial }{\partial z^{i}}+\overline{\xi }^{i+}%
\frac{\partial }{\partial \overline{z}^{i}})+e^{-}(\xi ^{i-}\frac{\partial }{%
\partial z^{i}}+\overline{\xi }^{i-}\frac{\partial }{\partial \overline{z}%
^{i}}); \\ 
z^{i}=\,z^{i+}e^{+}+z^{i-}e^{-};\,\overset{.}{\,z}^{i}=\overset{.}{\,z}%
^{i+}e^{+}+\overset{.}{z}^{i-}e^{-}=\xi ^{i+}e^{+}+\xi ^{i-}e^{-}; \\ 
\overline{z}^{i}=\,\overline{z}^{i+}e^{+}+\overline{z}^{i-}e^{-};\,\overset{.%
}{\,\overline{z}}^{i}=\overset{.}{\,\overline{z}}^{i+}e^{+}+\overset{.}{%
\overline{z}}^{i-}e^{-}=\overline{\xi }^{i+}e^{+}+\overline{\xi }^{i-}e^{-};%
\end{array}
\label{3.1}
\end{equation}%
where the dot indicates the derivative with respect to time $t$. The vector
field denoted by $V=(W^{+}-W^{-})(\xi )$ and given by%
\begin{equation}
(W^{+}-W^{-})(\xi )=e^{+}(-e^{\lambda }\xi ^{i+}\frac{\partial }{\partial
z^{i}}+e^{-\lambda }\overline{\xi }^{i+}\frac{\partial }{\partial \overline{z%
}^{i}})-e^{-}(-e^{\lambda }\xi ^{i-}\frac{\partial }{\partial z^{i}}%
+e^{-\lambda }\overline{\xi }^{i-}\frac{\partial }{\partial \overline{z}^{i}}%
)  \label{3.2}
\end{equation}%
is called conformal \textit{bi}-\textit{para} \textit{Liouville vector field}
on the bi-Lagrangian conformal manifold. The maps given by $T,$ $%
P:M\rightarrow A$ such that $T=\frac{1}{2}m_{i}(\overline{z}^{i})^{2}=\frac{1%
}{2}m_{i}(\overset{.}{z}^{i})^{2},$ $P=m_{i}gh$ are called \textit{the
kinetic energy} and \textit{the potential energy of the system,}
respectively.\textit{\ }Here\textit{\ }$m_{i},g$ and $h$ stand for mass of a
mechanical system having $n_{i}$ particle, the gravity acceleration and
distance to the origin of a mechanical system on the bi-Lagrangian conformal
manifold,

respectively. Then $L:M\rightarrow A$ is a map that satisfies the conditions;

\textbf{i)} $L=T-P$ is a conformal \textit{bi-para} \textit{Lagrangian
function,}

\textbf{ii)} the function given by $E_{L}=V(L)-L$ is \textit{a conformal
bi-para energy function}.

The operator $i_{(W^{+}-W^{-})}$ induced by $W^{+}-W^{-}$ and shown by%
\begin{equation}
i_{W^{+}-W^{-}}\omega (Z_{1},Z_{2},...,Z_{r})=\sum_{i=1}^{r}\omega
(Z_{1},...,(W^{+}-W^{-})Z_{i},...,Z_{r})  \label{3.3}
\end{equation}%
is said to be \textit{vertical derivation, }where $\omega \in \wedge
^{r}{}M, $ $Z_{i}\in \chi (M).$ The \textit{vertical differentiation} $%
d_{(P^{+}-P^{-})}$ is defined by%
\begin{equation}
d_{(W^{+}-W^{-})}=[i_{(W^{+}-W^{-})},d]=i_{(W^{+}-W^{-})}d-di_{(W^{+}-W^{-})}
\label{3.4}
\end{equation}%
where $d$ is the usual exterior derivation. For an almost para-complex
structure $W^{+}-W^{-}$, the closed para-conformal K\"{a}hlerian form is the
closed 2-form given by $\Phi _{L}=-dd_{_{(W^{+}-W^{-})}}L$ such that%
\begin{equation}
d_{\left( W^{+}-W^{-}\right) }L=\mathbf{e}^{+}\left( -e^{\lambda }\frac{%
\partial L}{\partial \overline{z}^{i}}dz^{i}+e^{-\lambda }\frac{\partial L}{%
\partial z^{i}}d\overline{z}^{i}\right) -\mathbf{e}^{-}\left( -e^{\lambda }%
\frac{\partial L}{\partial \overline{z}^{i}}dz^{i}+e^{-\lambda }\frac{%
\partial L}{\partial z^{i}}d\overline{z}^{i}\right) :\mathcal{F}%
(M)\rightarrow \wedge ^{1}{}M  \label{3.5}
\end{equation}%
Let $\xi $ be the second order differential equations given by equation (\ref%
{3.1}) and%
\begin{equation}
\begin{array}{l}
i_{\xi }\Phi _{L}=\Phi _{L}(\xi ) \\ 
=-\mathbf{e}^{+}\xi ^{i+}e^{\lambda }\frac{\partial \lambda }{\partial z^{j}}%
\frac{\partial L}{\partial \overline{z}^{i}}dz^{j}+\mathbf{e}^{+}\xi
^{i+}e^{\lambda }\frac{\partial \lambda }{\partial z^{j}}\frac{\partial L}{%
\partial \overline{z}^{i}}dz^{i}-\mathbf{e}^{+}\xi ^{i+}e^{\lambda }\frac{%
\partial ^{2}L}{\partial z^{j}\partial \overline{z}^{i}}dz^{j}+\mathbf{e}%
^{+}\xi ^{i+}e^{\lambda }\frac{\partial ^{2}L}{\partial z^{j}\partial 
\overline{z}^{i}}dz^{i} \\ 
-\mathbf{e}^{+}\overline{\xi }^{i+}e^{-\lambda }\frac{\partial \lambda }{%
\partial z^{j}}\frac{\partial L}{\partial z^{i}}dz^{i}+\mathbf{e}^{+}%
\overline{\xi }^{i+}e^{-\lambda }\frac{\partial ^{2}L}{\partial
z^{j}\partial z^{i}}dz^{i}+\mathbf{e}^{+}\overline{\xi }^{i+}e^{\lambda }%
\frac{\partial \lambda }{\partial \overline{z}^{j}}\frac{\partial L}{%
\partial \overline{z}^{i}}dz^{i}+\mathbf{e}^{+}\overline{\xi }%
^{i+}e^{\lambda }\frac{\partial ^{2}L}{\partial \overline{z}^{j}\partial 
\overline{z}^{i}}dz^{i} \\ 
+\mathbf{e}^{-}\xi ^{i-}e^{\lambda }\frac{\partial \lambda }{\partial z^{j}}%
\frac{\partial L}{\partial \overline{z}^{i}}dz^{j}-\mathbf{e}^{-}\xi
^{i-}e^{\lambda }\frac{\partial \lambda }{\partial z^{j}}\frac{\partial L}{%
\partial \overline{z}^{i}}dz^{i}+\mathbf{e}^{-}\xi ^{i-}e^{\lambda }\frac{%
\partial ^{2}L}{\partial z^{j}\partial \overline{z}^{i}}dz^{j}-\mathbf{e}%
^{-}\xi ^{i-}e^{\lambda }\frac{\partial ^{2}L}{\partial z^{j}\partial 
\overline{z}^{i}}dz^{i} \\ 
+\mathbf{e}^{-}\overline{\xi }^{i-}e^{-\lambda }\frac{\partial \lambda }{%
\partial z^{j}}\frac{\partial L}{\partial z^{i}}dz^{i}-\mathbf{e}^{-}%
\overline{\xi }^{i-}e^{-\lambda }\frac{\partial ^{2}L}{\partial
z^{j}\partial z^{i}}dz^{i}-\mathbf{e}^{-}\overline{\xi }^{i-}e^{\lambda }%
\frac{\partial \lambda }{\partial \overline{z}^{j}}\frac{\partial L}{%
\partial \overline{z}^{i}}dz^{i}-\mathbf{e}^{-}\overline{\xi }%
^{i-}e^{\lambda }\frac{\partial ^{2}L}{\partial \overline{z}^{j}\partial 
\overline{z}^{i}}dz^{i} \\ 
\mathbf{e}^{+}\xi ^{i+}e^{-\lambda }\frac{\partial \lambda }{\partial z^{j}}%
\frac{\partial L}{\partial z^{i}}d\overline{z}^{i}-\mathbf{e}^{+}\xi
^{i+}e^{-\lambda }\frac{\partial ^{2}L}{\partial z^{j}\partial z^{i}}d%
\overline{z}^{i}-\mathbf{e}^{+}\xi ^{i+}e^{\lambda }\frac{\partial \lambda }{%
\partial \overline{z}^{j}}\frac{\partial L}{\partial \overline{z}^{i}}d%
\overline{z}^{j}-\mathbf{e}^{+}\xi ^{i+}e^{\lambda }\frac{\partial ^{2}L}{%
\partial \overline{z}^{j}\partial \overline{z}^{i}}d\overline{z}^{j} \\ 
-\mathbf{e}^{+}\overline{\xi }^{i+}e^{-\lambda }\frac{\partial \lambda }{%
\partial \overline{z}^{j}}\frac{\partial L}{\partial z^{i}}d\overline{z}^{i}+%
\mathbf{e}^{+}\overline{\xi }^{i+}e^{-\lambda }\frac{\partial \lambda }{%
\partial \overline{z}^{j}}\frac{\partial L}{\partial z^{i}}d\overline{z}^{j}+%
\mathbf{e}^{+}\overline{\xi }^{i+}e^{-\lambda }\frac{\partial ^{2}L}{%
\partial \overline{z}^{j}\partial z^{i}}d\overline{z}^{i}-\mathbf{e}^{+}%
\overline{\xi }^{i+}e^{-\lambda }\frac{\partial ^{2}L}{\partial \overline{z}%
^{j}\partial z^{i}}d\overline{z}^{j} \\ 
-\mathbf{e}^{-}\xi ^{i-}e^{-\lambda }\frac{\partial \lambda }{\partial z^{j}}%
\frac{\partial L}{\partial z^{i}}d\overline{z}^{i}+\mathbf{e}^{-}\xi
^{i-}e^{-\lambda }\frac{\partial ^{2}L}{\partial z^{j}\partial z^{i}}d%
\overline{z}^{i}+\mathbf{e}^{-}\xi ^{i-}e^{\lambda }\frac{\partial \lambda }{%
\partial \overline{z}^{j}}\frac{\partial L}{\partial \overline{z}^{i}}d%
\overline{z}^{j}+\mathbf{e}^{-}\xi ^{i-}e^{\lambda }\frac{\partial ^{2}L}{%
\partial \overline{z}^{j}\partial \overline{z}^{i}}d\overline{z}^{j} \\ 
+\mathbf{e}^{-}\overline{\xi }^{i-}e^{-\lambda }\frac{\partial \lambda }{%
\partial \overline{z}^{j}}\frac{\partial L}{\partial z^{i}}d\overline{z}^{j}-%
\mathbf{e}^{-}\overline{\xi }^{i-}e^{-\lambda }\frac{\partial \lambda }{%
\partial \overline{z}^{j}}\frac{\partial L}{\partial z^{i}}d\overline{z}^{j}-%
\mathbf{e}^{-}\overline{\xi }^{i-}e^{-\lambda }\frac{\partial ^{2}L}{%
\partial \overline{z}^{j}\partial z^{i}}d\overline{z}^{j}+\mathbf{e}^{-}%
\overline{\xi }^{i-}e^{-\lambda }\frac{\partial ^{2}L}{\partial \overline{z}%
^{j}\partial z^{i}}d\overline{z}^{j}.%
\end{array}
\label{3.7}
\end{equation}%
Since the closed conformal para K\"{a}hlerian form $\Phi _{L}$ on $M$ is in
a para-symplectic\emph{\ }structure, it is found that%
\begin{equation}
\begin{array}{l}
E_{L}=\mathbf{e}^{+}\left( -\xi ^{i+}e^{\lambda }\frac{\partial L}{\partial 
\overline{z}^{i}}+\overline{\xi }^{i+}e^{-\lambda }\frac{\partial L}{%
\partial z^{i}}\right) -\mathbf{e}^{-}\left( -\xi ^{i-}e^{\lambda }\frac{%
\partial L}{\partial \overline{z}^{i}}+\overline{\xi }^{i-}e^{-\lambda }%
\frac{\partial L}{\partial z^{i}}\right) -L%
\end{array}
\label{3.8}
\end{equation}%
and thus%
\begin{equation}
\begin{array}{l}
dE_{L}=\mathbf{e}^{+}\left[ -\xi ^{i+}e^{\lambda }\frac{\partial \lambda }{%
\partial z^{j}}\frac{\partial L}{\partial \overline{z}^{i}}dz^{j}-\xi
^{i+}e^{\lambda }\frac{\partial ^{2}L}{\partial z^{j}\overline{z}^{i}}dz^{j}-%
\overline{\xi }^{i+}e^{-\lambda }\frac{\partial \lambda }{\partial z^{j}}%
\frac{\partial L}{\partial z^{i}}dz^{j}+\overline{\xi }^{i+}e^{-\lambda }%
\frac{\partial ^{2}L}{\partial z^{j}\partial z^{i}}dz^{j}\right] \\ 
-\mathbf{e}^{-}\left[ -\xi ^{i-}e^{\lambda }\frac{\partial \lambda }{%
\partial z^{j}}\frac{\partial L}{\partial \overline{z}^{i}}dz^{j}-\xi
^{i-}e^{\lambda }\frac{\partial ^{2}L}{\partial z^{j}\partial \overline{z}%
^{i}}dz^{j}-\overline{\xi }^{i-}e^{-\lambda }\frac{\partial \lambda }{%
\partial z^{j}}\frac{\partial L}{\partial z^{i}}dz^{j}+\overline{\xi }%
^{i-}e^{-\lambda }\frac{\partial ^{2}L}{\partial z^{j}\partial z^{i}}dz^{j}%
\right] -\frac{\partial L}{\partial z^{j}}dz^{j} \\ 
+\mathbf{e}^{+}\left[ -\xi ^{i+}e^{\lambda }\frac{\partial \lambda }{%
\partial \overline{z}^{j}}\frac{\partial L}{\partial \overline{z}^{i}}d%
\overline{z}^{j}-\xi ^{i+}e^{\lambda }\frac{\partial ^{2}L}{\partial 
\overline{z}^{j}\partial \overline{z}^{i}}d\overline{z}^{j}-\overline{\xi }%
^{i+}e^{-\lambda }\frac{\partial \lambda }{\partial \overline{z}^{j}}\frac{%
\partial L}{\partial z^{i}}d\overline{z}^{j}+\overline{\xi }^{i+}e^{-\lambda
}\frac{\partial ^{2}L}{\partial \overline{z}^{j}\partial z^{i}}d\overline{z}%
^{j}\right] \\ 
-\mathbf{e}^{-}\left[ -\xi ^{i-}e^{\lambda }\frac{\partial \lambda }{%
\partial \overline{z}^{j}}\frac{\partial L}{\partial \overline{z}^{i}}d%
\overline{z}^{j}-\xi ^{i-}e^{\lambda }\frac{\partial ^{2}L}{\partial 
\overline{z}^{j}\partial \overline{z}^{i}}d\overline{z}^{j}-\overline{\xi }%
^{i-}e^{-\lambda }\frac{\partial \lambda }{\partial \overline{z}^{j}}\frac{%
\partial L}{\partial z^{i}}d\overline{z}^{j}+\overline{\xi }^{i-}e^{-\lambda
}\frac{\partial ^{2}L}{\partial \overline{z}^{j}\partial z^{i}}d\overline{z}%
^{j}\right] -\frac{\partial L}{\partial \overline{z}^{i}}d\overline{z}^{j}.%
\end{array}
\label{3.9}
\end{equation}%
Use of equation (\ref{1}) gives:%
\begin{equation}
\begin{array}{l}
\mathbf{e}^{+}\xi ^{i+}e^{\lambda }\frac{\partial ^{2}L}{\partial
z^{j}\partial \overline{z}^{i}}+\mathbf{e}^{+}\overline{\xi }^{i+}e^{\lambda
}\frac{\partial ^{2}L}{\partial \overline{z}^{j}\partial \overline{z}^{i}}-%
\mathbf{e}^{-}\xi ^{i-}e^{\lambda }\frac{\partial ^{2}L}{\partial
z^{j}\partial \overline{z}^{i}}-\mathbf{e}^{-}\overline{\xi }^{i-}e^{\lambda
}\frac{\partial ^{2}L}{\partial \overline{z}^{j}\partial \overline{z}^{i}}
\\ 
+\mathbf{e}^{+}\xi ^{i+}e^{\lambda }\frac{\partial \lambda }{\partial z^{j}}%
\frac{\partial L}{\partial \overline{z}^{i}}+\mathbf{e}^{+}\overline{\xi }%
^{i+}e^{\lambda }\frac{\partial \lambda }{\partial \overline{z}^{j}}\frac{%
\partial L}{\partial \overline{z}^{i}}-\mathbf{e}^{-}\xi ^{i-}e^{\lambda }%
\frac{\partial \lambda }{\partial z^{j}}\frac{\partial L}{\partial \overline{%
z}^{i}}-\mathbf{e}^{-}\overline{\xi }^{i-}e^{\lambda }\frac{\partial \lambda 
}{\partial \overline{z}^{j}}\frac{\partial L}{\partial \overline{z}^{i}}+%
\frac{\partial L}{\partial z^{j}} \\ 
-\mathbf{e}^{+}\xi ^{i+}e^{-\lambda }\frac{\partial ^{2}L}{\partial
z^{j}\partial z^{i}}-\mathbf{e}^{+}\overline{\xi }^{i+}e^{-\lambda }\frac{%
\partial ^{2}L}{\partial \overline{z}^{j}\partial z^{i}}+\mathbf{e}^{-}\xi
^{i-}e^{-\lambda }\frac{\partial ^{2}L}{\partial z^{j}\partial z^{i}}+%
\mathbf{e}^{-}\overline{\xi }^{i-}e^{-\lambda }\frac{\partial ^{2}L}{%
\partial \overline{z}^{j}\partial z^{i}} \\ 
\mathbf{e}^{+}\xi ^{i+}e^{-\lambda }\frac{\partial \lambda }{\partial z^{j}}%
\frac{\partial L}{\partial z^{i}}+e^{+}\overline{\xi }^{i+}e^{-\lambda }%
\frac{\partial \lambda }{\partial \overline{z}^{j}}\frac{\partial L}{%
\partial z^{i}}-\mathbf{e}^{-}\xi ^{i-}e^{-\lambda }\frac{\partial \lambda }{%
\partial z^{j}}\frac{\partial L}{\partial z^{i}}-e^{-}\overline{\xi }%
^{i-}e^{-\lambda }\frac{\partial \lambda }{\partial \overline{z}^{j}}\frac{%
\partial L}{\partial z^{i}}+\frac{\partial L}{\partial \overline{z}^{j}}=0.%
\end{array}
\label{3.10}
\end{equation}%
If a curve denoted by $\alpha :A\rightarrow M$ is considered to be an
integral curve of $\xi ,$ $\xi (L)=\frac{\partial L}{\partial t}$ $,$ then
the following equation is obtained:%
\begin{equation}
\begin{array}{l}
\left( \mathbf{e}^{+}-\mathbf{e}^{-}\right) e^{\lambda }\left[ \mathbf{e}%
^{+}\xi ^{i+}\frac{\partial }{\partial z^{j}}+\mathbf{e}^{+}\overline{\xi }%
^{i+}\frac{\partial }{\partial \overline{z}^{j}}+\mathbf{e}^{-}\xi ^{i-}%
\frac{\partial }{\partial z^{j}}+\mathbf{e}^{-}\overline{\xi }^{i-}\frac{%
\partial }{\partial \overline{z}^{j}}\right] \left( \frac{\partial L}{%
\partial \overline{z}^{j}}\right) \\ 
+\left( \mathbf{e}^{+}-\mathbf{e}^{-}\right) e^{\lambda }\left[ \mathbf{e}%
^{+}\xi ^{i+}\frac{\partial }{\partial z^{j}}+\mathbf{e}^{+}\overline{\xi }%
^{i+}\frac{\partial }{\partial \overline{z}^{j}}+\mathbf{e}^{-}\xi ^{i-}%
\frac{\partial }{\partial z^{j}}+\mathbf{e}^{-}\overline{\xi }^{i-}\frac{%
\partial }{\partial \overline{z}^{j}}\right] \left( \lambda \right) \left( 
\frac{\partial L}{\partial \overline{z}^{j}}\right) +\frac{\partial L}{%
\partial z^{j}} \\ 
-\left( \mathbf{e}^{+}-\mathbf{e}^{-}\right) e^{-\lambda }\left[ \mathbf{e}%
^{+}\xi ^{i+}\frac{\partial }{\partial z^{j}}+\mathbf{e}^{+}\overline{\xi }%
^{i+}\frac{\partial }{\partial \overline{z}^{j}}+\mathbf{e}^{-}\xi ^{i-}%
\frac{\partial }{\partial z^{j}}+\mathbf{e}^{-}\overline{\xi }^{i-}\frac{%
\partial }{\partial \overline{z}^{j}}\right] \frac{\partial L}{\partial z^{j}%
} \\ 
\left( \mathbf{e}^{+}-\mathbf{e}^{-}\right) e^{-\lambda }\left[ \mathbf{e}%
^{+}\xi ^{i+}\frac{\partial }{\partial z^{j}}+\mathbf{e}^{+}\overline{\xi }%
^{i+}\frac{\partial }{\partial \overline{z}^{j}}+\mathbf{e}^{-}\xi ^{i-}%
\frac{\partial }{\partial z^{j}}+\mathbf{e}^{-}\overline{\xi }^{i-}\frac{%
\partial }{\partial \overline{z}^{j}}\right] \left( \lambda \right) \frac{%
\partial L}{\partial z^{j}}+\frac{\partial L}{\partial \overline{z}^{j}}=0%
\end{array}
\label{3.11}
\end{equation}%
or%
\begin{eqnarray}
\left( \mathbf{e}^{+}-\mathbf{e}^{-}\right) e^{\lambda }\xi \left( \frac{%
\partial L}{\partial \overline{z}^{i}}\right) +\left( \mathbf{e}^{+}-\mathbf{%
e}^{-}\right) e^{\lambda }\xi \left( \lambda \right) \left( \frac{\partial L%
}{\partial \overline{z}^{i}}\right) +\frac{\partial L}{\partial z^{i}} &=&0,
\label{3.12} \\
-\left( \mathbf{e}^{+}-\mathbf{e}^{-}\right) e^{-\lambda }\xi \left( \frac{%
\partial L}{\partial z^{i}}\right) +\left( \mathbf{e}^{+}-\mathbf{e}%
^{-}\right) e^{-\lambda }\xi \left( \lambda \right) \frac{\partial L}{%
\partial z^{i}}+\frac{\partial L}{\partial \overline{z}^{i}} &=&0.  \notag
\end{eqnarray}%
Then the following equations are found:%
\begin{equation}
\left( \mathbf{e}^{+}-\mathbf{e}^{-}\right) \frac{\partial }{\partial t}%
\left( e^{\lambda }\frac{\partial L}{\partial \overline{z}^{i}}\right) +%
\frac{\partial L}{\partial z^{i}}=0\text{ \ , \ }\left( \mathbf{e}^{+}-%
\mathbf{e}^{-}\right) \frac{\partial }{\partial t}\left( e^{-\lambda }\frac{%
\partial L}{\partial z^{i}}\right) -\frac{\partial L}{\partial \overline{z}%
^{i}}=0.  \label{3.13}
\end{equation}%
Thus equations\textbf{\ }(\ref{3.13}) are seen to be \textbf{conformal
bi-para Euler-Lagrange equations} on the distributions $D_{1}$ and $D_{2},$
and then the triple $(M,\Phi _{L},\xi )$ is seen to be a \textbf{conformal
bi-para mechanical system}\textit{\ }with taking into account the basis $%
\{e^{+},e^{-}\}$ on the bi-Lagrangian conformal manifold $(M,\Phi
,D_{1},D_{2})$.

\section{Conformal Bi-Para Hamiltonians}

In the part, conformal bi-para Hamilton equations and a conformal bi-para
Hamiltonian mechanical system on the bi Lagrangian conformal manifold $%
(M,\Phi ,D_{1},D_{2})$ will be derived. Let $(z_{i},\overline{z}_{i})$ be
its para-complex coordinates. Let $\{\frac{\partial }{\partial z_{i}},\frac{%
\partial }{\partial \overline{z}_{i}}\}$ and $\{dz_{i},d\overline{z}_{i}\},$
be bases and cobases of $T_{p}(M)$ and $T_{p}^{\ast }(M)$ of $M,$
respectively. Let us assume that an almost bi-para-complex conformal
structure, a bi-para-conformal Liouville form and a bi-para-complex
conformal 1-form on the distributions ${}D_{1}$ and $D_{2}$ are shown by $%
W^{\ast +}-W^{\ast -}$, $\lambda $ and $\omega $, respectively. Then, we
using \cite{miron} and (\ref{2.11}):

\begin{equation}
\begin{array}{c}
\omega =\frac{1}{2}[(z_{i}dz_{i}+\overline{z}_{i}d\overline{z}%
_{i})e^{+}+(e^{2\lambda }z_{i}dz_{i}+e^{2\lambda }\overline{z}_{i}d\overline{%
z}_{i})e^{-}], \\ 
\lambda =(W^{\ast +}-W^{\ast -})(\omega )=\frac{1}{2}[(\mathbf{-}%
e^{+}e^{\lambda }z_{i}d\overline{z}_{i}+e^{+}e^{\lambda }\overline{z}%
_{i}dz_{i})] \\ 
-\frac{1}{2}(\mathbf{-}e^{-}e^{\lambda }z_{i}d\overline{z}%
_{i}+e^{-}e^{\lambda }\overline{z}_{i}dz_{i})].%
\end{array}
\label{4.111}
\end{equation}%
It is well known that if $\Phi $ is a closed para-K\"{a}hlerian form on the
bi-Lagrangian conformal manifold, then $\Phi $ is also a para-symplectic
structure on ${}$the bi-Lagrangian conformal manifold. Given a
bi-para-conformal Hamiltonian vector field $Z_{H}$ fixed with the
bi-para-conformal Hamiltonian energy\textit{\ }$H$ that is 
\begin{equation}
Z_{H}=(Z_{i}\frac{\partial }{\partial z_{i}}+\overline{Z}_{i}\frac{\partial 
}{\partial \overline{z}_{i}})e^{+}+(Z_{i}\frac{\partial }{\partial z_{i}}+%
\overline{Z}_{i}\frac{\partial }{\partial \overline{z}_{i}})e^{-}.
\label{4.2}
\end{equation}%
Then closed 2-form is%
\begin{eqnarray}
\Phi &=&-d\lambda =e^{+}\overline{Z}_{i}dz_{i}-e^{-}\overline{Z}_{i}dz_{i}-%
\frac{1}{2}\left[ 
\begin{array}{c}
e^{+}e^{\lambda }\frac{\partial \lambda }{\partial z_{i}}\bar{Z}%
_{i}z_{i}dz_{i}+e^{+}e^{\lambda }\frac{\partial \lambda }{\partial \overline{%
z}_{i}}\overline{Z}_{i}\overline{z}_{i}dz_{i} \\ 
-e^{-}e^{\lambda }\frac{\partial \lambda }{\partial z_{i}}\overline{Z}%
_{i}z_{i}dz_{i}-e^{-}e^{\lambda }\frac{\partial \lambda }{\partial \overline{%
z}_{i}}\overline{Z}_{i}\overline{z}_{i}dz_{i}%
\end{array}%
\right]  \label{4.3} \\
&&-e^{+}Z_{i}d\overline{z}_{i}+e^{-}Z_{i}d\overline{z}_{i}-\frac{1}{2}\left[ 
\begin{array}{c}
\mathbf{-}e^{+}e^{\lambda }\frac{\partial \lambda }{\partial z_{i}}%
Z_{i}z_{i}d\overline{z}_{i}-e^{+}e^{\lambda }\frac{\partial \lambda }{%
\partial \overline{z}_{i}}Z_{i}\overline{z}_{i}d\overline{z}_{i} \\ 
+e^{-}e^{-}e^{\lambda }\frac{\partial \lambda }{\partial z_{i}}Z_{i}z_{i}d%
\overline{z}_{i}+e^{-}e^{\lambda }\frac{\partial \lambda }{\partial 
\overline{z}_{i}}Z_{i}\overline{z}_{i}d\overline{z}_{i}%
\end{array}%
\right] .  \notag
\end{eqnarray}%
And then it follows%
\begin{equation}
\begin{array}{c}
i_{Z_{H}}\Phi =\Phi (Z_{H}) \\ 
=\bar{Z}_{i}e^{+}\left[ 1-\frac{1}{2}e^{\lambda }\left( z_{i}\frac{\partial
\lambda }{\partial z_{i}}+\overline{z}_{i}\frac{\partial \lambda }{\partial 
\overline{z}_{i}}\right) \right] dz_{i}-Z_{i}e^{+}\left[ 1+\frac{1}{2}%
e^{\lambda }\left( z_{i}\frac{\partial \lambda }{\partial z_{i}}+\overline{z}%
_{i}\frac{\partial \lambda }{\partial \overline{z}_{i}}\right) \right] d%
\overline{z}_{i} \\ 
-\bar{Z}_{i}e^{-}\left[ 1-\frac{1}{2}e^{\lambda }\left( z_{i}\frac{\partial
\lambda }{\partial z_{i}}+\overline{z}_{i}\frac{\partial \lambda }{\partial 
\overline{z}_{i}}\right) \right] dz_{i}+Z_{i}e^{-}\left[ 1+\frac{1}{2}%
e^{\lambda }\left( z_{i}\frac{\partial \lambda }{\partial z_{i}}+\overline{z}%
_{i}\frac{\partial \lambda }{\partial \overline{z}_{i}}\right) \right] d%
\overline{z}_{i}%
\end{array}
\label{4.4}
\end{equation}%
On the other hand, the differential of the bi-para-conformal Hamiltonian
energy $H$\textit{\ }is calculated as follows:%
\begin{equation}
dH=(\frac{\partial H}{\partial z_{i}}dz_{i}+\frac{\partial H}{\partial 
\overline{z}_{i}}d\overline{z}_{i})e^{+}+(\frac{\partial H}{\partial z_{i}}%
dz_{i}+\frac{\partial H}{\partial \overline{z}_{i}}d\overline{z}_{i})e^{-}.
\label{4.5}
\end{equation}%
By means of equation\textbf{\ }(\ref{3}), using equation (\ref{4.4}) and (%
\ref{4.5}), the conformal bi-para Hamiltonian vector field is seen to be%
\begin{equation}
Z_{H}=(Z_{i}\frac{\partial }{\partial z_{i}}+\overline{Z}_{i}\frac{\partial 
}{\partial \overline{z}_{i}})e^{+}+(Z_{i}\frac{\partial }{\partial z_{i}}+%
\overline{Z}_{i}\frac{\partial }{\partial \overline{z}_{i}})e^{-}.
\label{4.6}
\end{equation}%
If a curve $\alpha :I\subset A\rightarrow M$ is an integral curve of the
conformal bi-para Hamiltonian vector field $Z_{H}$, i.e., $Z_{H}(\alpha (t))=%
\dot{\alpha}(t)$ $,\,\,t\in I.$ In the local coordinates, we get $\alpha
(t)=(z_{i}(t),\overline{z}_{i}(t))$ and%
\begin{equation}
\dot{\alpha}(t)=(\frac{dz_{i}}{dt}\frac{\partial }{\partial z_{i}}+\frac{d%
\overline{z}_{i}}{dt}\frac{\partial }{\partial \overline{z}_{i}})e^{+}+(%
\frac{dz_{i}}{dt}\frac{\partial }{\partial z_{i}}+\frac{d\overline{z}_{i}}{dt%
}\frac{\partial }{\partial \overline{z}_{i}})e^{-}.  \label{4.10}
\end{equation}%
Taking equations\textbf{\ }(\ref{4.6}),\textbf{\ }(\ref{4.10}), the
following equations are found%
\begin{equation}
\frac{dz_{i}}{dt}=\frac{-(e^{+}-e^{-})}{\left[ 1+\frac{1}{2}e^{\lambda
}\left( z_{i}\frac{\partial \lambda }{\partial z_{i}}+\overline{z}_{i}\frac{%
\partial \lambda }{\partial \overline{z}_{i}}\right) \right] }\frac{\partial
H}{\partial \overline{z}_{i}}\text{ , }\frac{d\overline{z}_{i}}{dt}=\frac{%
(e^{+}-e^{-})}{\left[ 1-\frac{1}{2}e^{\lambda }\left( z_{i}\frac{\partial
\lambda }{\partial z_{i}}+\overline{z}_{i}\frac{\partial \lambda }{\partial 
\overline{z}_{i}}\right) \right] }\frac{\partial H}{\partial z_{i}}.
\label{4.12}
\end{equation}%
Hence, equations\textbf{\ }(\ref{4.12}) are seen to be \textbf{conformal
bi-para Hamilton equations} on the bi-Lagrangian conformal manifold $(M,\Phi
,D_{1},D_{2}),$ and then the triple $(M,\Phi ,Z_{H})$ is seen to be a 
\textbf{conformal bi-para Hamiltonian mechanical system}\textit{\ }with the
use of basis $\{e^{+},e^{-}\}$ on $(M,\Phi ,D_{1},D_{2})$.

\section{Conclusion}

It is seen in the above, formalisms of Lagrangian and Hamiltonian mechanics
had intrinsically been described by taking into account the basis $%
\{e^{+},e^{-}\}$\ on the bi-Lagrangian conformal manifold $(M,\Phi
,D_{1},D_{2})$. Conformal bi-para Lagrangian and bi-para Hamiltonian models
have arisen to be very important tools since they present a simple method to
describe the model for bi-para-conformal mechanical systems. So, the
equations derived in (\ref{3.13}) and (\ref{4.12}) are only considered to be
a first step to realize how bi-para-complex conformal geometry has been used
in solving problems in different physically spaces. For further research,
bi-para-complex conformal Lagrangian and Hamiltonian vector fields derived
here are suggested to deal with problems in different fields of physics. In
the literature, the equations, which explains the linear orbits of the
objects, were presented. This study explained the non-linear orbits of the
objects in the space by the help of revised equations using Weyl theorem.

\end{document}